\documentclass[10pt]{amsart}
\usepackage{amssymb}

\newtheorem{thm1}{Theorem}[section]
\newtheorem{lem1}[thm1]{Lemma}

\newtheorem{def1}[thm1]{Definition}

\newtheorem{prop1}[thm1]{Proposition}
\newtheorem{ex1}[thm1]{Example}

\begin{document}

\title[Radical of a lattice ideal]
{Binomial generation of the radical of a lattice ideal}
\author[A. Katsabekis]{Anargyros Katsabekis}
\address {Department of Mathematics, Aristotle University of Thessaloniki, Thessaloniki 54124, Greece} \email{katsabek@aegean.gr}
\author[M. Morales]{ Marcel
Morales}
\address { Universit\'e de Grenoble I, Institut Fourier,
UMR 5582, B.P.74, 38402 Saint-Martin D'H\`eres Cedex, and IUFM de
Lyon, 5 rue Anselme, 69317 Lyon Cedex France}
\email{morales@ujf-grenoble.fr}
\author[A. Thoma]{ Apostolos Thoma}
\address { Department of Mathematics, University of Ioannina,
Ioannina 45110, GREECE } \email{athoma@cc.uoi.gr}


\begin{abstract}
\par Let $I_{L, \rho}$ be a lattice ideal. We provide a necessary and sufficient criterion under which a set of
binomials in $I_{L, \rho}$ generate the radical of $I_{L, \rho}$
up to radical. We apply our results to the problem of determining
the minimal number of generators of $I_{L, \rho}$ or of the
$rad(I_{L, \rho})$ up to radical.
\end{abstract}
\maketitle

\section{Introduction}
\par Lattice ideals is an important class of binomial ideals  with a
lot of applications in several areas like algebraic statistics,
dynamical systems, graph theory, hypergeometric differential
equations, integer programming and toric geometry \cite{E-S}.
Lattice ideals are  generalizations of toric ideals, see \cite{MS,
St} for details about these ideals. The generation of the radical
of a lattice ideal up to radical by binomials was the subject of
several recent papers, see \cite{B, BL, BMTh, BMTh2, BJT, BGT,
E-S, EV, Eto, Eto1, Eto2, KMT, KT, MT} and sections of books
\cite{V1, V2}. In \cite{EV} S. Eliahou and R. Villarreal, see
Theorem 2.5, provide two necessary and sufficient conditions for a
set of binomials to generate a toric ideal up to radical. This
result was later generalized by K. Eto \cite{Eto, Eto1} for
lattice ideals. The last years appeared results in the literature
\cite{B, KMT} that required a huge number of binomials  to
generate the radical of lattice ideals up to radical. In this
article we approach the problem in a different manner to
understand these results.
\par Let $(L,\rho )$ be a partial character on $\mathbb{Z}^m$, then we associate to any lattice ideal $I_{L, \rho}$ a rational polyhedral cone $\sigma_{L}=pos_{\mathbb{Q}}(A)$, where $A=\{{\bf a}_{1},\ldots,{\bf a}_{m}\}$, see section 2 for details. Let $S$ be a subset of the cone
$\sigma_L$, then ${\mathbb{E}}_S:=\{i \in \{1, \dots ,m\}|{\bf
a}_i\in  S\}$. Given a vector ${\bf u}=(u_1,\ldots,u_m) \in
\mathbb{Z}^m$, its {\em support} is the set $${\rm supp}({\bf
u})=\{i \in \{1,\ldots,m\} \ | \ u_{i} \neq 0\}.$$ We can write
every vector ${\bf u}$ in $\mathbb{Z}^{m}$ uniquely as ${\bf
u}={\bf u}_{+}-{\bf u}_{-}$, where ${\bf u}_{+}$ and ${\bf u}_{-}$
are non-negative and have disjoint support. We shall denote by
$\mathbb{Z}^E$ the set $\{{\bf u} \in \mathbb{Z}^m \ | \ {\rm
supp}({\bf u}) \subset E\}$.
\begin{def1} {\rm We say that the set of binomials $\{B({\bf u}_1), B({\bf u}_2),
\dots , B({\bf u}_q)\} \subset I_{L,\rho }$ is a cover of $A$ if
and only if for every $E\subset \{1,\dots ,m\}$, that is not in
the form $\mathbb{E}_{\mathcal{F}}$ for a face $\mathcal{F}$ of
$\sigma_L$, there exist an $i \in \{1,\ldots,q\}$ such that $({\bf
u}_i)_+ \in \mathbb{Z}^E$ and $({\bf u}_i)_- \notin \mathbb{Z}^E$
or $({\bf u}_i)_- \in \mathbb{Z}^E$ and $({\bf u}_i)_+ \notin
\mathbb{Z}^E$.}
\end{def1}
The next theorem improves and generalizes results of S. Eliahou-R.
Villarreal \cite{EV} and K. Eto \cite{Eto, Eto1}. It shows how the
binomial generation of the radical of a lattice ideal
 is related with the geometry of the cone  $\sigma_L$ and the algebra of the lattice $L$.
\bigskip
\newline
\textbf{Theorem 3.5.} Let $(L,\rho)$ be a partial character on
$\mathbb{Z}^m$, ${\bf u}_{1},\ldots,{\bf u}_{q}$ be elements of
the lattice $L$ and $\sigma_L={pos}_{\mathbb{Q}}(A)$ the rational
polyhedral cone associated to $I_{L,\rho }$. Then
$$rad(I_{L,\rho })=rad(B({\bf u}_1), B({\bf u}_2), \dots ,
B({\bf u}_q))$$ if and only if \begin{enumerate} \item[(i)]
$\{B({\bf u}_1), B({\bf u}_2), \dots , B({\bf u}_q)\}$ is a cover
of $A$, \item[(ii)] for every face $\mathcal{F}$ of $\sigma_L $ we
have
$$L \cap
\mathbb{Z}^{\mathbb{E}_{\mathcal{F}}}=\sum _{{\bf u}_i \in
\mathbb{Z}^{\mathbb{E}_{\mathcal{F}}}} \mathbb{Z}{\bf u}_i,$$ in
characteristic zero. While in characteristic $p>0$,
$$(L \cap
\mathbb{Z}^{\mathbb{E}_{\mathcal{F}}}):p^{\infty}=(\sum _{{\bf
u}_i\in \mathbb{Z}^{\mathbb{E}_{\mathcal{F}}}}\mathbb{Z}{\bf
u}_i):p^{\infty}.$$
\end{enumerate}
The two conditions depend also on the geometry of the cone
$\sigma_{L}$. In particular our first condition depend on the non
faces of $\sigma_{L}$ and the second involve sublattices
associated to faces of the cone. The new conditions explain why
there are lattices that require a huge number of binomials to
generate the radical of a lattice ideal up to radical,  see the
remark after Proposition 3.6. In section 4 we provide a necessary
and sufficient condition for a lattice ideal to be complete
intersection. Finally we prove that lattice ideals with associated
full vector configuration are either set-theoretic complete
intersection or almost set-theoretic complete intersection in
characteristic zero, while they are set-theoretic complete
intersection in positive characteristic.\\

\section{Basics on lattice ideals}

\par Let $K$ be a field  of any characteristic. A {\em lattice} is a finitely generated free abelian group. A
{\em partial character} $(L,\rho )$ on $\mathbb{Z}^m$ is a
homomorphism $\rho $ from a sublattice $L$ of $\mathbb{Z}^m$ to
the multiplicative group $K^*=K-\{ 0 \}$. Given a partial
character $(L, \rho )$ on $\mathbb{Z}^m$, we define the ideal
$$I_{L,\rho }:=(\{B({\bf \alpha }):={\bf x}^{{\bf \alpha }_+}-
\rho({\bf \alpha } ){\bf x}^{{\bf \alpha } _-}|{\bf \alpha } ={\bf
\alpha } _+-{\bf \alpha} _- \in L\})\subset K[x_1,\dots ,x_m]$$
called {\em lattice ideal}. Where ${\bf x}^{{\bf \beta
}}=x_1^{\beta_1}\cdots x_m^{\beta_m}$ for ${\bf \beta}
=(\beta_1,\dots ,\beta_m)\in \mathbb{N}^m$. The height of
$I_{L,\rho }$ is equal to the rank of $L$. Given a finite subset
$C$ of $L$, define the ideal
$$J_{{C},\rho}:=(\{B({\bf \alpha }) \ | \ {\bf \alpha } ={\bf
\alpha } _+-{\bf \alpha} _- \in C\}) \subset I_{L,\rho }.$$
\begin{lem1}\cite {St} A subset $C$ spans the lattice $L$
if and only if $$J_{{C},\rho}:(x_1\cdots
x_m)^{\infty}=I_{L,\rho}.$$
\end{lem1}

If $L$ is a sublattice of $\mathbb{Z}^m$, then the {\em
saturation} of $L$ is the lattice
$$ Sat(L):=\{  {\bf \alpha } \in \mathbb{Z}^m |d {\bf \alpha } \in L \
\textrm{for some} \ d \in \mathbb{Z}, \ d \neq 0 \}.$$ We say that
the lattice $L$ is {\em saturated} if $L=Sat(L)$. The lattice
ideal $I_{L,\rho }$
is prime if and only if $L$ is saturated.  \\
For a prime number $p$ the $p$-saturation of $L$ is the lattice
$$(L:p^{\infty}):=\{{\bf \alpha } \in \mathbb{Z}^m |p^k {\bf \alpha } \in L \
\textrm{for some} \ k \in \mathbb{N}\}.$$
\par Throughout this paper we assume that $L$ is a non-zero  {\em positive sublattice}
 of $\mathbb{Z}^m$, that is $L\cap \mathbb{N}^m=\{{\bf 0}\}$. This means that the
lattice ideal $I_{L,\rho }$ is
homogeneous with respect to some positive grading.\\
The  group $\mathbb{Z}^m/Sat(L)$ is free abelian, therefore is
isomorphic to $\mathbb{Z}^n$,
 where $n=m-rank(L)$. Let $\psi $ be the above isomorphism, ${\bf e}_1,\dots ,{\bf e}_m$
the unit vectors of $\mathbb{Z}^m$ and ${\bf a}_i:=\psi ({\bf
e}_i+Sat(L))$ for $1\leq i\leq m$. We call $A=\{{\bf a}_i| 1\leq
i\leq
m \} \subset \mathbb{Z}^n$ {\em the configuration of vectors associated to the lattice} $L$. \\
We grade $K[x_1,\dots ,x_m]$ by setting ${\rm deg}_A(x_i)={\bf
a}_i$ for $i=1,\dots ,m$. We define the {\em $A$-degree} of the
monomial ${\bf x^u}$  to be
$${\rm deg}_A ({\bf x^u}) := u_1{\bf a}_1+\cdots +u_{m}{\bf
a}_m\in \mathbb{N} A,$$ where $\mathbb{N}A $ is the semigroup
generated by $A$. A polynomial $F\in K[x_1,\dots ,x_m]$ is called
$A$-{\em homogeneous} if the monomials in each non zero term of
$F$ have the same $A$-degree. An ideal $I$ is $A$-homogeneous if
it is generated by $A$-homogeneous polynomials. The lattice ideal
$I_{L, \rho }$, as well as every binomial in it, is
$A$-homogeneous.

The {\em binomial arithmetical rank} of a binomial ideal $I$
(written ${\rm bar}(I)$) is the smallest integer $s$ for which
there exist binomials $f_1,\dots ,f_s$ in $I$ such that $rad(I) =
rad(f_1,\ldots ,f_s)$. Hence the binomial arithmetical rank is an
upper bound for the {\em arithmetical rank} of a binomial ideal
(written ${\rm ara}(I)$), which is the smallest integer $s$ for
which there exists $f_1,\ldots ,f_s$ in $I$ such that $rad(I) =
rad(f_1,\ldots ,f_s)$. From the definitions, the generalized
Krull's principal ideal theorem and the graded version of
Nakayama's Lemma we deduce the following inequality for a lattice
ideal $I_{L,\rho }$:
$$h(I_{L, \rho})\leq {\rm
ara}(I_{L, \rho}) \leq {\rm bar}(I_{L, \rho})\leq \mu(I_{L,
\rho}).$$ Here $h(I)$ denotes the height and $\mu(I)$ denotes the
minimal number of generators of an ideal $I$. When $h(I) = {\rm
ara}(I)$ the ideal $I$ is called a {\em set-theoretic complete
intersection}; when $h(I) = \mu(I)$ it is called a {\em complete
intersection}. The ideal $I$ is called an {\em almost
set-theoretic complete intersection} if ${\rm ara}(I) \leq
h(I)+1$.

We associate to the lattice ideal $I_{L,\rho }$ the rational
polyhedral cone $$\sigma_L =pos_{\mathbb{Q}}(A):=\{\lambda_1{\bf
a}_1+\cdots +\lambda_m{\bf a}_m| \lambda_i\in \mathbb{Q} \
\textrm{and} \ \lambda_i \geq 0\}.$$ The {\em dimension} of
$\sigma_L$ is equal to the dimension of the $\mathbb{Q}$-vector
space $${\rm span}_{\mathbb{Q}}(\sigma_L)=\{\lambda_1{\bf
a}_1+\cdots +\lambda_m{\bf a}_m| \lambda_i\in \mathbb{Q}\}.$$ A
{\em face} of $\sigma_L$ is any set of the form
$$\mathcal{F}=\sigma_L \cap \{{\bf x} \in \mathbb{Q}^n: {\bf c}{\bf
x}=0\}$$ where ${\bf c} \in \mathbb{Q}^{n}$ and ${\bf c}{\bf x}
\geq 0$ for all $x \in \sigma_L$. Faces of dimension one are
called {\em extreme rays}. If the number of the extreme rays of a
cone coincides with the dimension (i.e. the extreme rays are
linearly independent), the cone is called {\em simplex cone}. A
cone $\sigma$ is {\em strongly convex} if $\sigma \cap -\sigma
=\{{\bf 0}\}$. The condition that the lattice $L$ is positive, is
equivalent with the condition that the cone $\sigma_L$ is strongly
convex.

We decompose the affine space $K^m$ into $2^m$ {\em coordinate
cells}, $$(K^*)^{E}:=\{(q_1, \ldots ,q_m)\in K^m | q_i \not= 0 \
for \ i\in {E}, q_i=0\ for \ i\notin {E} \},$$ where ${E}$ runs
over all subsets of $\{ 1,\ldots ,m \}$. Let $P=(x_1,\dots ,x_m)$
be a point of $K^m$ then
$$P_{E}:=(\delta^{{E}}_1x_1,
\delta^{{E}}_2x_2, \ldots , \delta^{{E}}_mx_m)\in K^m, $$ where
$\delta^{{E}}_i=1$ if $i\in E$ and $\delta^{{E}}_i=0$ if $i\notin
{E}$. Note that if $P\in (K^*)^{ \{ 1,\ldots ,m \} }$ then
$P_{E}\in (K^*)^{E}$.

The $n$-dimensional algebraic torus $( K^*)^n$ acts on the affine
$m$-space $ K^m$ via $$(x_1, \ldots ,x_m) \rightarrow (x_1{\bf
t}^{{\bf a}_1}, \ldots ,x_m{\bf t}^{{\bf a}_m}).$$ Let
$\mathcal{K}$ denotes the algebraic closure of $K$. For the
lattice algebraic set $V(I_{L, \rho})\subset {\mathcal{K}}^m$ we
have that $V(I_{L,\rho })=\cup _{j=1}^{g}{\bf X}_{A,j}$, see
\cite{E-S} and \cite{KMT}, where the affine toric variety ${\bf
X}_{A,j}$ is the Zariski-closure of the $( \mathcal{K}^*)^n$-orbit
of a point $P_{j}=(c_{j1},c_{j2},\ldots ,c_{jm})\in \mathcal{K}^m$
for appropriate $c_{ji}$ all different from zero. Actually the
toric variety ${\bf X}_{A,j}$ is the disjoint union of the orbits
of the points $(P_j)_{\mathbb{E}_{\mathcal F}}$, for every face
${\mathcal F}$ of $\sigma_L$. There are no points of the toric
varieties ${\bf X}_{A,j}$ which are in the cells
$(\mathcal{K}^*)^{E}$, where $E$ is not in the form
$\mathbb{E}_{\mathcal F}$ for a face ${\mathcal F}$ of $\sigma_L$.
Thus the lattice algebraic set $V(I_{L, \rho})$ has points only on
the cells in the form $(\mathcal{K}^*)^{\mathbb{E}_{\mathcal{F}}}$
for some face $\mathcal{F}$ of the cone $\sigma_L$.\\

\section{Radical generation by binomials}

Set $K[E]=K[\{x_{i} | i \in E\}]$, where $E$ is a subset of
$\{1,\ldots,m\}$.

\begin{lem1} {\rm (Cf. \cite{KMT}, Lemma 2.3)} Let $\mathcal{F}$ be
a face of the rational polyhedral cone
$\sigma_L=pos_{\mathbb{Q}}(A)$. The monomial ${\bf x}^{\bf u} \in
K[\mathbb{E}_{\mathcal{F}}]$ if and only if ${\rm deg}_{A}({\bf
x}^{\bf u}) \in \mathcal{F}$.
\end{lem1}
To every face $\mathcal{F}$ of $\sigma_L$ we can associate the
ideal $I_{L_{\mathcal{F}},\rho}:=I_{L,\rho } \cap
 K[\mathbb{E}_{\mathcal{F}}]$. The next proposition shows how the
 geometry of the cone affects the radical generation by
 A-homogeneous polynomials. Note that binomials in a lattice ideal
 are always A-homogeneous.
\begin{prop1}Let $B=\{G_{1}, G_{2}, \ldots ,G_{q}\} \subset I_{L,\rho }$
be a set of $A$-homogenous polynomials and
$\sigma_L=pos_{\mathbb{Q}}(A)$ the rational polyhedral cone
associated to $I_{L, \rho}$. Then $B$ generate $rad(I_{L, \rho})$
up to radical if and only if for every face $\mathcal{F}$ of
$\sigma_L$ the set $B\cap K[\mathbb{E}_{\mathcal{F}}]$ generate
$rad(I_{L_{\mathcal{F},\rho}})$ up to radical.
\end{prop1}

\noindent \textbf{Proof.} ($\Leftarrow$) For
$\mathcal{F}=\sigma_L$ we have that
$K[\mathbb{E}_{\sigma_L}]=K[x_{1},\ldots,x_{m}]$, so
$I_{L_{\sigma_{L}},\rho}=I_{L,\rho}$ and also $B \cap
K[\mathbb{E}_{\sigma_{L}}]=B$.\\
($\Rightarrow$) Obviously $B\cap K[\mathbb{E}_{\mathcal{F}}]
\subset I_{L_{\mathcal{F},\rho}}$ and therefore $rad(B\cap
K[\mathbb{E}_{\mathcal{F}}]) \subset
rad(I_{L_{\mathcal{F}},\rho})$. It is enough to show that
$I_{L_{\mathcal{F}},\rho} \subset rad(B \cap
K[\mathbb{E}_{\mathcal{F}}])$. Let $f \in I_{L_{\mathcal{F}},\rho}
\subset I_{L,\rho}$, so there exist a positive integer $k$ and
polynomials $A_1, \ldots A_q \in K[x_1,\ldots, x_m]$ such that
$$f^k=A_1G_{1}+\cdots +A_qG_{q}.$$ For every
polynomial $h \in K[x_1,\ldots,x_m]$, we let $h_{\mathcal{F}}$ the
polynomial in $K[\mathbb{E}_{\mathcal{F}}]$ taken from $h$ by
setting all the variables $x_i$, where $i \notin
\mathbb{E}_{\mathcal{F}}$, equal to zero. There are two cases for
the polynomial $G_{i}$:
\begin{enumerate} \item $G_{i}$ belongs to
$K[\mathbb{E}_{\mathcal{F}}]$. In this case
$(G_{i})_{\mathcal{F}}=G_{i}$. \item the polynomial $G_{i}$ does
not belong to $K[\mathbb{E}_{\mathcal{F}}]$. Then, from Lemma 3.1,
every monomial of $G_{i}$ does not belong to
$K[\mathbb{E}_{\mathcal{F}}]$ since the polynomial $G_{i}$ is
$A$-homogeneous. In this case $(G_{i})_{\mathcal{F}}=0$.
\end{enumerate}
It follows that $B\cap K[\mathbb{E}_{\mathcal{F}}]$ consists of
the nonzero polynomials among
$(G_{1})_{\mathcal{F}},\ldots,(G_{q})_{\mathcal{F}}$. We have that
$$f^k=(f_{\mathcal{F}})^k=(A_{1})_{\mathcal{F}}(G_{1})_{\mathcal{F}}+\cdots +(A_{q})_{\mathcal{F}} (G_{q})_{\mathcal{F}},$$ since $f \in
I_{L_{\mathcal{F}},\rho}=I_{L,\rho} \cap
K[\mathbb{E}_{\mathcal{F}}]$. Thus $f^{k}$ belongs to the ideal
generated by the set $B\cap K[\mathbb{E}_{\mathcal{F}}]$. \hfill
$\square$

\begin{lem1} Let $(L, \rho)$ be a partial character on $\mathbb{Z}^m$ and $\mathcal{F}$ a face of $\sigma_L=
pos_{\mathbb{Q}}(A)$. Then a vector ${\bf u}={\bf u}_+ -{\bf u}_-
\in L$ belongs to $\mathbb{Z}^{\mathbb{E}_{\mathcal{F}}}$ if and
only if ${\bf u}_{+} \in \mathbb{Z}^{\mathbb{E}_{\mathcal{F}}}$ if
and only if ${\bf u}_{-} \in
\mathbb{Z}^{\mathbb{E}_{\mathcal{F}}}$.
\end{lem1}
\noindent \textbf{Proof.} It suffices to prove ${\bf u}_{+} \in
\mathbb{Z}^{\mathbb{E}_{\mathcal{F}}}$ implies that ${\bf u}_{-}
\in \mathbb{Z}^{\mathbb{E}_{\mathcal{F}}}$. The other parts
are similar or obvious.\\
From Lemma 3.1 we have that ${\rm deg}_{A}({\bf x}^{{\bf u}_+})$
belongs to $\mathcal{F}$ and therefore ${\rm deg}_{A}({\bf
x}^{{\bf u}_-})$ belongs to $\mathcal{F}$, since the binomial
$B({\bf u})$ is $A$-homogeneous.\\ Let ${\bf
u}_-=(u_{-,1},\ldots,u_{-,m})$. For any vector ${\bf
c}_{\mathcal{F}}$ which defines the face $\mathcal{F}$ we have
that ${\bf c}_{\mathcal{F}} \cdot (\sum_{i=1}^{m} u_{-,i} {\bf
a}_i)=0$. Thus $\sum_{i=1}^{m}u_{-,i} ({\bf c}_{\mathcal{F}} \cdot
{\bf a}_i)=0$. But $u_{-,i} \geq 0$ for all $i$, and ${\bf
c}_{\mathcal{F}} \cdot {\bf a}_i>0$ for every $i \not \in
\mathbb{E}_{\mathcal{F}}$, while ${\bf c}_{\mathcal{F}} \cdot {\bf
a}_i=0$ for every $i\in \mathbb{E}_{\mathcal{F}}$. So we conclude
that $u_{-,i}=0$ for every $i \not\in \mathbb{E}_{\mathcal{F}}$.
Consequently ${\bf u}_{-} \in
\mathbb{Z}^{\mathbb{E}_{\mathcal{F}}}$. \hfill $\square$

\begin{prop1} Let $(L, \rho)$ be a partial character on $\mathbb{Z}^m$. For a face $\mathcal{F}$ of
$\sigma_L=pos_{\mathbb{Q}}(A)$, the ideals
$I_{L_{\mathcal{F}},\rho}$ and $I_{{L\cap
\mathbb{Z}^{\mathbb{E}_{\mathcal{F}}}},\rho}$ coincide.
\end{prop1}
\noindent \textbf{Proof.} From Corollary 1.3 in \cite{E-S} the
ideal $I_{L_{\mathcal{F}},\rho}$ has a generating set $$\{B({\bf
u}_1),\ldots,B({\bf u}_q)\} \subset I_{L,\rho } \cap
K[\mathbb{E}_{\mathcal{F}}].$$ So ${\bf u}_{i}$ belongs to $L$,
for every $i=1,\ldots,q$. If ${\bf u}_{i} \notin
\mathbb{Z}^{\mathbb{E}_{\mathcal{F}}}$, then, from Lemma 3.3, both
$({\bf u}_i)_{+} \notin \mathbb{Z}^{\mathbb{E}_{\mathcal{F}}}$,
$({\bf u}_i)_{-} \notin \mathbb{Z}^{\mathbb{E}_{\mathcal{F}}}$ and
therefore $B({\bf u}_{i})$ is not in $I_{L_{\mathcal{F}},\rho}$, a
contradiction, so ${\bf u}_i \in L\cap
\mathbb{Z}^{\mathbb{E}_{\mathcal{F}}}$. Conversely if $B({\bf v})$
belongs to $I_{{L\cap
\mathbb{Z}^{\mathbb{E}_{\mathcal{F}}}},\rho}$, then ${\bf v}$
belongs to $L$ and therefore $B({\bf v}) \in I_{L}$. Moreover
${\bf v} \in \mathbb{Z}^{\mathbb{E}_{\mathcal{F}}}$, so, again
from Lemma 3.3, both ${\bf v}_{+}$ and ${\bf v}_{-}$ belong to
$\mathbb{Z}^{\mathbb{E}_{\mathcal{F}}}$. Thus $B({\bf v}) \in
K[\mathbb{E}_{\mathcal{F}}]$ and therefore belongs to
$I_{L_{\mathcal{F}},\rho}$. \hfill $\square$\\

The next theorem generalizes results of S. Eliahou-R. Villarreal
\cite{EV} and K. Eto \cite{Eto, Eto1}.
\begin{thm1} Let $(L,\rho)$ be a partial character on
$\mathbb{Z}^m$, ${\bf u}_{1},\ldots,{\bf u}_{q}$ be elements of
the lattice $L$ and $\sigma_L={pos}_{\mathbb{Q}}(A)$ the rational
polyhedral cone associated to $I_{L,\rho }$. Then
$$rad(I_{L,\rho })=rad(B({\bf u}_1), B({\bf u}_2), \dots ,
B({\bf u}_q))$$ if and only if \begin{enumerate} \item[(i)]
$\{B({\bf u}_1), B({\bf u}_2), \dots , B({\bf u}_q)\}$ is a cover
of $A$, \item[(ii)] for every face $\mathcal{F}$ of $\sigma_L $ we
have
$$L \cap
\mathbb{Z}^{\mathbb{E}_{\mathcal{F}}}=\sum _{{\bf u}_i \in
\mathbb{Z}^{\mathbb{E}_{\mathcal{F}}}} \mathbb{Z}{\bf u}_i,$$ in
characteristic zero. While in characteristic $p>0$,
$$(L \cap
\mathbb{Z}^{\mathbb{E}_{\mathcal{F}}}):p^{\infty}=(\sum _{{\bf
u}_i\in \mathbb{Z}^{\mathbb{E}_{\mathcal{F}}}}\mathbb{Z}{\bf
u}_i):p^{\infty}.$$
\end{enumerate}
\end{thm1}
\noindent \textbf{Proof.} ($\Rightarrow$) (i) Let $E\subset
\{1,\dots ,m\}$ that is not in the form $\mathbb{E}_{\mathcal{F}}$
for any face $\mathcal{F}$ of $\sigma_L$. Then for all $j$ the
point $(P_j)_E$ is not a point of $V(I_{L,\rho })$, since
$(P_j)_E$ belongs to the cell $(\mathcal{K}^*)^{E}$ and $E$ is of
the above form. Thus $(P_j)_E$ is not a point of $V(B({\bf u}_1),
B({\bf u}_2), \dots , B({\bf u}_q))$ and therefore there is an $i
\in \{1,\ldots,q\}$ such that $B({\bf u}_{i})((P_j)_E) \neq 0$. On
the other hand $$B({\bf u}_{i})(P_j)=0, \ \textrm{since} \ P_j \in
{\bf X}_{A,j}.$$ If ${\rm supp}({\bf u}_i)\subset E$ then we have
$$B({\bf u}_{i})((P_j)_E)=B({\bf u}_{i})(P_j)=0,$$ a contradiction. So ${\rm supp}({\bf
u}_i)\not\subset E$ which implies that
$${\bf x}^{({\bf u}_i)_-}((P_j)_E)=0 \ \ \textrm{or} \ \ {\bf x}^{({\bf
u}_i)_+}((P_j)_E)=0 .$$ Since $B({\bf u}_{i})((P_j)_E) \neq 0$ we
have
$${\bf x}^{({\bf u}_i)_+}((P_j)_E)\not=0 \ \ \textrm{and} \ \ {\bf
x}^{({\bf u}_i)_-}((P_j)_E)=0 \ \ \textrm{or}$$ $${\bf x}^{({\bf
u}_i)_+}((P_j)_E)=0 \ \ \textrm{and} \ \ {\bf x}^{({\bf
u}_i)_-}((P_j)_E)\not= 0.$$ Thus $({\bf u}_i)_+ \in \mathbb{Z}^E$
and $({\bf u}_i)_- \not\in \mathbb{Z}^E$ or $({\bf u}_i)_- \in
\mathbb{Z}^E$ and $({\bf u}_i)_+ \not\in \mathbb{Z}^E$.
Consequently $\{B({\bf
u}_1), B({\bf u}_2), \dots , B({\bf u}_q)\}$ is a cover of $A$.\\
(ii) Let $\mathcal{F}$ be a face of $\sigma_L$. Set $G=\sum _{{\bf
u}_i \in \mathbb{Z}^{\mathbb{E}_{\mathcal{F}}}} \mathbb{Z}{\bf
u}_i$, $B=\{B({\bf u}_1),\ldots, B({\bf u}_q)\}$ and $B \cap
K[\mathbb{E}_{\mathcal{F}}]=\{B({\bf u}_{i_1}),\ldots,B({\bf
u}_{i_t})\}$. In any characteristic we have $I_{G, \rho} \subset
I_{{L\cap \mathbb{Z}^{\mathbb{E}_{\mathcal{F}}}},\rho}$ and also
$$(B({\bf u}_{i_1}),\ldots,B({\bf
u}_{i_t})) \subset I_{G, \rho} \subset I_{{L \cap
\mathbb{Z}^{\mathbb{E}_{\mathcal{F}}}},\rho}.$$ Thus
$$rad(B({\bf u}_{i_1}),\ldots,B({\bf
u}_{i_t})) \subset rad(I_{G, \rho}) \subset rad(I_{{L\cap
\mathbb{Z}^{\mathbb{E}_{\mathcal{F}}}},\rho}).$$But $rad(I_{L,\rho
})=rad(B({\bf u}_1), B({\bf u}_2), \dots , B({\bf u}_q))$, so
combining Propositions 3.2 and 3.4, since every binomial in
$I_{L,\rho }$ is $A$-homogeneous, we deduce that $$rad(I_{{L\cap
\mathbb{Z}^{\mathbb{E}_{\mathcal{F}}}},\rho})=rad(B({\bf
u}_{i_1}),\ldots,B({\bf u}_{i_t}))$$ and therefore $$rad(I_{G,
\rho})=rad(B({\bf u}_{i_1}),\ldots,B({\bf u}_{i_t})).$$ Thus
$rad(I_{{L\cap
\mathbb{Z}^{\mathbb{E}_{\mathcal{F}}}},\rho})=rad(I_{G, \rho})$.
In characteristic zero we have that $I_{G, \rho}=I_{{L\cap
\mathbb{Z}^{\mathbb{E}_{\mathcal{F}}}},\rho}$, so $L\cap
\mathbb{Z}^{\mathbb{E}_{\mathcal{F}}}=G$, while in characteristic
$p>0$ it holds $I_{G:p^{\infty}, \rho}=I_{{L\cap
\mathbb{Z}^{\mathbb{E}_{\mathcal{F}}}:p^{\infty}},\rho}$, so
$G:p^{\infty}=L\cap
\mathbb{Z}^{\mathbb{E}_{\mathcal{F}}}:p^{\infty}$, see \cite{E-S}.\\
($\Leftarrow$) It is enough to prove that $V(B({\bf u}_1), B({\bf
u}_2), \dots , B({\bf u}_q)) \subset V(I_{L, \rho})$. Let ${\bf y}
\in V(B({\bf u}_1), B({\bf u}_2), \ldots , B({\bf u}_q)) \subset
\mathcal{K}^{m}$ and assume that ${\bf y} \in
(\mathcal{K}^*)^{E}$, where $E$ is not of the form
$\mathbb{E}_{\mathcal{F}}$ for a face $\mathcal{F}$ of $\sigma_L$.
Then there is an $i \in \{1,\ldots,q\}$ such that $({\bf u}_i)_+
\in \mathbb{Z}^E$ and $({\bf u}_i)_- \not\in \mathbb{Z}^E$ since
$\{B({\bf u}_1), B({\bf u}_2), \dots , B({\bf u}_q)\}$ is a cover
of $A$. This means that $B({\bf u}_{i})({\bf y}) \neq 0$, a
contradiction to the fact that ${\bf y} \in V(B({\bf u}_1), B({\bf
u}_2), \dots , B({\bf u}_q))$. Thus ${\bf y} \in
(\mathcal{K}^*)^{\mathbb{E}_{\mathcal{F}}}$, for a face
$\mathcal{F}$ of $\sigma_L$. Let ${\bf v} \in L$. We will prove
that $B({\bf v})({\bf y})=0$. If ${\bf v} \notin
\mathbb{Z}^{\mathbb{E}_{\mathcal{F}}}$, then by Lemma 3.3 both
${\bf v}_{+}$, ${\bf v}_{-}$ do not belong to
$\mathbb{Z}^{\mathbb{E}_{\mathcal{F}}}$ and therefore $B({\bf
v})({\bf y})=0$. Suppose now that ${\bf v} \in
\mathbb{Z}^{\mathbb{E}_{\mathcal{F}}}$, which implies that $B({\bf
v}) \in I_{{L\cap \mathbb{Z}^{\mathbb{E}_{\mathcal{F}}}}, \rho}$.
Set $C=\{{\bf u}_i | {\rm supp}({\bf u}_i) \subset
\mathbb{E}_{\mathcal{F}}\} \subset
\mathbb{Z}^{\mathbb{E}_{\mathcal{F}}}$. Then
\begin{enumerate} \item by hypothesis in characteristic zero we
have that $C$ spans the lattice $L\cap
\mathbb{Z}^{\mathbb{E}_{\mathcal{F}}}$. Thus, from Lemma 2.1, the
binomial $B({\bf v})$ belongs to the ideal $J_{{C},\rho}:<\prod_{j
\in \mathbb{E}_{\mathcal{F}}}x_{j}>^{\infty}$. But ${\bf y}$
belongs to $V(B({\bf u}_1), B({\bf u}_2),\ldots ,B({\bf u}_q))$,
which in particular implies that ${\bf y} \in V(B({\bf u}_i)| {\bf
u}_i \in C)$. Consequently $B({\bf v})({\bf y})=0$. \item by
hypothesis in characteristic $p>0$ we have that
$$(L\cap
\mathbb{Z}^{\mathbb{E}_{\mathcal{F}}}):p^{\infty}=(\sum_{{\bf u}_i
\in C}\mathbb{Z}{\bf u}_i):p^{\infty},$$ so $I_{({L\cap
\mathbb{Z}^{\mathbb{E}_{\mathcal{F}}}):p^{\infty}},\rho}=I_{G:p^{\infty},
\rho}$ and therefore $rad(I_{{L\cap
\mathbb{Z}^{\mathbb{E}_{\mathcal{F}}}},\rho})=rad(I_{G, \rho})$
where $G=\sum_{{\bf u}_i \in C}\mathbb{Z}{\bf u}_i$. Thus $(B({\bf
v}))^{k} \in I_{G, \rho}$, so $B({\bf v})(y)=0$ since $I_{G,
\rho}=J_{{C},\rho}:p^{\infty}$. \hfill $\square$\\
\end{enumerate}

It is well known, see Proposition 8.7 in \cite{E-S} and
Proposition 3.4 in \cite{Eto}, that every lattice ideal is the
set-theoretic intersection of its circuits. Therefore a cover of
$A$ always exists. In the case that the $\sigma_L$ is a simplex
cone it is easy to find a cover. For an ${\bf a}_{i}\in A$, $1
\leq i \leq m$, we define ${\mathcal{ F}}_{{\bf a}_i}$ to be the
{\em minimal face} of $\sigma_L $ that contains ${\bf a}_i$, i.e.
$${\mathcal F}_{{\bf a}_i}=\cap_{{\bf a}_i \in {\mathcal F}}{\mathcal
F},$$ since any intersection of faces of $\sigma_L$ is a face of
$\sigma_L$.
\begin{prop1} Let $(L,\rho)$ be a partial character on
$\mathbb{Z}^m$ such that the rational polyhedral cone
$\sigma_L=pos_{\mathbb{Q}}(A) \subset \mathbb{Q}^n$ associated to
$I_{L,\rho }$ is a simplex cone of dimension $n$, then there exist
a cover of $A$ consisting of $m-n$ binomials.
\end{prop1}
\noindent \textbf{Proof.} The strongly convex cone $\sigma_L$ is
$n$-dimensional and also it is a simplex cone, so we choose one
vector for each extreme ray of $\sigma_L $ and obtain a linearly
independent set $B \subset A$ consisting of $n$ vectors such that
$\sigma_L=pos_{\mathbb{Q}}(B)$. Remark that the choice of the
above vectors is not unique. Let $A=\{{\bf a}_1,\ldots,{\bf
a}_m\}$. We rearrange the vectors in $A$ such that the first $n$
vectors, i.e. ${\bf a}_1,\ldots,{\bf a}_n$, are the elements of
$B$. Every ${\mathcal{ F}}_{{\bf a}_i}$ is a rational polyhedral
cone, so ${\mathcal{ F}}_{{\bf a}_i}$ is of the form
$pos_{\mathbb{Q}}({\bf a}_j| j\in \mathbb{E}_{{ \mathcal{F}}_{{\bf
a}_i}} \cap \{1,\dots ,n\})$. For every $i \in \{n+1,\ldots,m\}$
the vector ${\bf a}_i$ belongs to ${\mathcal{ F}}_{{\bf a}_i}$, so
there exist a binomial $B({\bf u}_i)=x_i^{g_i}-N_i \in I_{L,\rho
}$ where $x_i^{g_i}={\bf x}^{({\bf u}_i)_{+}}$, $N_i={\bf
x}^{({\bf u}_i)_{-}}$ and the support of $({\bf u}_i)_{-}$ equals
$\mathbb{E}_{{ \mathcal{F}}_{{\bf a}_i}} \cap \{1,\dots ,n\}$. We
claim that the set $\{B({\bf u}_{n+1}),\dots ,B({\bf
u}_m)\}$ is a cover of $A$.\\
Let $E\subset \{1,\ldots,n,n+1,\ldots,m\}$ which is not of the
form $\mathbb{E}_{\mathcal{F}}$ for a face $\mathcal{F}$ of the
cone $\sigma_L $. Let $\mathcal{H}_E=pos_{\mathbb{Q}}({\bf
a}_j|j\in E\cap \{1,\dots ,n\})$, then $\mathcal{H}_E$ is a face
of the simplex cone $\sigma_L $. We have that $E\cap \{1,\dots
,n\} \subsetneqq E$ and therefore $E$ has at least one element
belonging to the set $\{n+1,\ldots,m\}$. There are two cases:
\begin{enumerate} \item The set $E$ has an element $i \in
\{n+1,\ldots,m\}$ such that ${\bf a}_i \notin \mathcal{H}_E$. Then
$\mathbb{E}_{{ \mathcal{F}}_{{\bf a}_i}} \cap \{1,\dots ,n\}$ is
not a subset of $E$, since if it was the vector ${\bf a}_i$ should
belong to $\mathcal{H}_E$. For the binomial $B({\bf u}_i)$ we have
 $({\bf u}_i)_+ \in \mathbb{Z}^E$ and $({\bf u}_i)_-
\notin \mathbb{Z}^E$. \item For every element $i \in
\{n+1,\ldots,m\}$ of $E$ the vector ${\bf a}_i \in \mathcal{H}_E$.
In this case there exist a $j \in \{n+1,\ldots,m\}$, which does
not belong to $E$, such that ${\bf a}_j \in \mathcal{H}_E$. If
not, we have that $\mathcal{H}_E=pos_{\mathbb{Q}}({\bf a}_j | j
\in E)$ and therefore $E=\mathbb{E}_{\mathcal{H}_E}$
contradiction. For the binomial $B({\bf u}_j) \in I_{L,\rho }$ we
have $({\bf u}_j)_- \in \mathbb{Z}^E$, since $\mathcal{F}_{{\bf
a}_j} \subset \mathcal{H}_E$, and $({\bf u}_j)_+ \notin
\mathbb{Z}^E$.
\end{enumerate} The preceding discussion yields that $\{B({\bf u}_{n+1}),\dots ,B({\bf u}_m)\}$ is a cover of $A$.
\hfill $\square$\\

Let us explain why there are lattices which require a huge number
of binomials to generate the radical of a lattice ideal up to
radical. The first condition of Theorem 3.5 states that you need a
number of binomials to cover $A$. In the case that $\sigma_L$ is a
simplex cone of dimension $n$ Proposition 3.6 provides a cover
consisting of $m-n$ binomials. But if the geometry of the cone is
complicated, then the number of binomials increases. In \cite{KMT}
we studied an example of a family of toric ideals of dimension $n$
and proved, by explicit computation of $5 {{n}\choose {3}}+6
{{n}\choose {4}}$ binomials generating the toric ideal up to
radical, that the binomial arithmetical rank is equal to $5
{{n}\choose {3}}+6 {{n}\choose {4}}$. The above binomials
constitute a cover and actually one can easily prove that there is
no cover with less than $5 {{n}\choose {3}}+6 {{n}\choose {4}}$
binomials. The second condition depends also on the characteristic
of the field and it states that you do not only need to generate
the lattice or the p-saturation of it in the characteristic $p$
case, but also sublattices or their p-saturations that are
associated to the faces of the cone $\sigma_L$. The examples given
by M. Barile in \cite{B} show this effect very clearly and the
dependance also on the characteristic.
 In \cite{B} M. Barile provide a family of
toric ideals of dimension $n$ and height $n$ that their cones
$\sigma _L$ are simplex cones, so to cover $A$ you need only $n$
binomials, see proposition 3.6. While to generate the lattice and
the sublattices associated to the faces of the cone $\sigma_L $
you need $n+ {{n-1}\choose {2}}$ binomials when the characteristic
of the field is not equal to $p$. In characteristic $p$ exactly
$n$ elements of the lattice provide a cover and generate also the
p-saturations of all the sublattices associated to the faces of
the cone $\sigma_L$. In this case the toric ideals are binomial
set theoretic complete intersections and the lattice is completely
$p$-glued \cite{MT}.

To compute the exact value of the binomial arithmetical rank is
not usually an easy process, since one can use the same or modify
some binomials to satisfy more than one conditions of Theorem 3.5,
as the following example shows. Note  that the procedure may
depend also on the characteristic of the field.

\begin{ex1} {\rm In this example we apply our methods to compute the exact
value of the binomial arithmetical rank of a lattice ideal. For a
different approach to the same example see also \cite{Eto1}. Given
an $n \times m$ matrix $M$ with columns the vectors of the set
$A=\{{\bf a}_1,\ldots,{\bf a}_{m}\} \subset \mathbb{Z}^n$ and a
saturated lattice $L=ker_{\mathbb{Z}}(M) \subset \mathbb{Z}^m$,
the toric ideal $I_{L,1}$ will be denoted by $I_A$. For a face
$\mathcal{F}$ of $pos_{\mathbb{Q}}(A)$ we shall denote by
$A_{\mathcal{F}}$ the set $\{{\bf a}_i | i \in
\mathbb{E}_{\mathcal{F}}\}$. We consider the set of vectors
$A=\{(a_{1},a_{2},a_{3}) \in \mathbb{N}^{3} \ | \
a_{1}+a_{2}+a_{3}=3\}$. The vectors of $A$ are the transpose of
the columns of the matrix $$M_{3,3}=\begin{pmatrix} 3 & 0 & 0 & 2
& 1 & 0 & 0 & 2 & 1 & 1
            \\ 0 & 3 & 0 & 1 & 2 & 2 & 1 & 0 & 0 & 1
            \\ 0 & 0 & 3 & 0 & 0 & 1 & 2 & 1 & 2 & 1
            \end{pmatrix}.$$
Consider the toric ideal $I_{A} \subset K[x_{1},\ldots,x_{10}]$.
The toric variety $V(I_A) \subset K^{10}$
is the so called $(3,3)$-Veronese toric variety.\\
The cone $pos_{\mathbb{Q}}(A)$ is a three dimensional simplicial
cone with three facets ${\mathcal{F}_1}, {\mathcal{F}_2},
{\mathcal{F}_3}$.
 For
$\mathbb{E}_{\mathcal{F}_1}=\{1,2,4,5\}$ we can see, using
\cite{marcel}, that the set $C_{1}=\{{\bf u}_1,{\bf u}_2,{\bf
u}_3\}$ spans the lattice $L \cap
\mathbb{Z}^{\mathbb{E}_{\mathcal{F}_1}}$ where ${\bf
u}_1=(2,1,0,-3,0,0,0,0,0,0)$, ${\bf u}_2=(1,0,0,-2,1,0,0,0,0,0)$
and ${\bf u}_{3}=(1,2,0,0,-3,0,0,0,0,0)$. Moreover
$$I_{A_{\mathcal{F}_1}}=(B({\bf u}_1),B({\bf u}_2),B({\bf u}_3)).$$
For $\mathbb{E}_{\mathcal{F}_2}=\{2,3,6,7\}$ we can see, using
\cite{marcel}, that the set $C_{2}=\{{\bf u}_4,{\bf u}_5,{\bf
u}_6\}$ spans the lattice $L \cap
\mathbb{Z}^{\mathbb{E}_{\mathcal{F}_2}}$ where ${\bf
u}_4=(0,2,1,0,0,-3,0,0,0,0)$, ${\bf u}_5=(0,1,0,0,\\0,-2,1,0,0,0)$
and ${\bf u}_6=(0,1,2,0,0,0,-3,0,0,0)$. In fact
$$I_{A_{\mathcal{F}_2}}=(B({\bf u}_4),B({\bf u}_5),B({\bf u}_6)).$$
For $\mathbb{E}_{\mathcal{F}_3}=\{1,3,8,9\}$ we use \cite{marcel}
to deduce that the set $C_{3}=\{{\bf u}_7,{\bf u}_8,{\bf u}_9\}$
spans the lattice $L \cap \mathbb{Z}^{\mathbb{E}_{\mathcal{F}_3}}$
where ${\bf u}_7=(2,0,1,0,0,0,0,-3,0,0)$, ${\bf
u}_8=(1,0,0,0,0,0,0,\\-2,1,0)$ and ${\bf
u}_9=(1,0,2,0,0,0,0,0,-3,0)$. In fact
$$I_{A_{\mathcal{F}_3}}=(B({\bf u}_7),B({\bf u}_8),B({\bf u}_9)).$$
For $\mathbb{E}_{\sigma}=\{1,\ldots,10\}$ we can see that $$L=L
\cap \mathbb{Z}^{\mathbb{E}_{\sigma}}=\mathbb{Z}{\bf
u}_1+\mathbb{Z}{\bf u}_2+\mathbb{Z}{\bf u}_4+\mathbb{Z}{\bf
u}_5+\mathbb{Z}{\bf u}_8+\mathbb{Z}{\bf u}_{10}+\mathbb{Z}{\bf
u}_{11}$$ where ${\bf u}_{10}=(1,0,0,-1,0,0,0,-1,0,1)$ and ${\bf
u}_{11}=(0,0,0,1,0,0,1,0,0,-2)$. For the singleton $E=\{10\}$
there is a binomial, namely $B({\bf u}_{11})$, in the toric ideal
$I_{A}$ such that $({\bf u}_{11})_{-} \in \mathbb{Z}^E$ and $({\bf
u}_{11})_{+} \notin \mathbb{Z}^E$. Therefore we can check that the
set $$\{B({\bf u}_{1}),B({\bf u}_{3}),B({\bf u}_{4}),B({\bf
u}_{6}),B({\bf u}_{7}),B({\bf u}_{9}),B({\bf u}_{11})\}$$ is a
cover of $A$. Thus
$rad(I_{A})=rad(B({\bf u}_{1}),\ldots,B({\bf u}_{11}))$  and so  ${\rm bar}(I_{A})\leq 11$.\\
Suppose that $K$ is a field of characteristic zero. We will prove
that ${\rm bar}(I_{A})=11$. Let $\mathcal{B}$ be a set of
binomials which generate $I_{A}$ up to radical. Using the fact
that every $I_{A_{\mathcal{F}_i}}$ is generated up to radical by
$3$ binomials and it is not a set-theoretic complete intersection
on binomials see \cite{BMTh2}, we take, from Proposition 3.2, that
$\mathcal{B}$ has at least $9$ binomials. Let $B({\bf v}_i^1),
B({\bf v}_i^2), B({\bf v}_i^3) \in I_{\mathcal{F}_i}$, $1 \leq i
\leq 3$, be those $9$ binomials. Note that in all these binomials
the variable $x_{10}$ does not appear.  Since the set
$\mathcal{B}$ must cover the set $\{10\}$ it must contain also a
monic binomial $B({\bf w}) \in I_{A}$ in $x_{10}$. Note that for
such ${\bf w}=(w_1,\ldots,w_{10}) \in \mathbb{Z}^{10}$ we have
$w_{10}>1$. Consider the vector ${\bf z}=(1,0,0,-1,0,0,0,-1,0,1)
\in L$. It does not belong to the lattice generated by the $9$
vectors ${\bf v}_i^j$ plus the vector ${\bf w}$, since the last
coordinate is $1$.  Consequently in $\mathcal{B}$ there are more
than $10$ binomials.
 Therefore ${\rm bar}(I_{A})\geq 11$, we conclude that ${\rm bar}(I_{A})= 11$. \\
In the case that $K$ is a field of characteristic $p>0$ we have
that ${\rm ara}(I_{A_{\mathcal{F}_i}})={\rm
bar}(I_{A_{\mathcal{F}_i}})=2$, see \cite{BMTh}. In fact when $K$
is a field of characteristic $3$ we have that ${\rm bar}(I_{A})=7$
since
$$rad(I_{A})=rad(B({\bf u}_1),B({\bf u}_3),B({\bf
u}_4),B({\bf u}_6),B({\bf u}_7),B({\bf u}_9),B({\bf v})),$$ where
${\bf v}=(1,1,1,0,0,0,0,0,0,-3)$, and the height of $I_{A}$ equals
$7$. This means also that the semigroup generated by $A$ is
completely 3-glued, see \cite{BMTh2}.
 When $K$
is a field of characteristic $p\not= 3$ the toric ideal $I_A$ is
generated up to radical by $8$ binomials.\\}
\end{ex1}

\section{Applications}

Complete intersection lattice ideals have been characterized in
\cite{FMS}, \cite{MT}, either in terms of semigroup gluing or in
terms of mixed dominating matrices. Both characterizations show
that the problem of determining the complete intersection property
for lattice ideals is in the NP-class \cite{FMS}, \cite{SSS}.
Therefore it is interesting to find better criteria for
establishing that a lattice ideal is or is not a complete
intersection. In this direction such criteria were given in
\cite{MiT} which can be read off from the geometry of the cone
$\sigma_L $. Our next Theorem provides a criterion depending on
the lattices associated to the faces of the cone $\sigma_L $. It
generalizes in a more geometric setting a result by K. Eto (see
Lemma 1.6 in \cite{Eto2})
 for complete intersections lattice ideals and the proof, contrary to the proof in \cite{Eto2}, uses mixed
dominating matrices. Note that lattice divisors correspond to
facets of the cone $\sigma_L $.
\begin{thm1} Let $L$ be a nonzero positive sublattice of
$\mathbb{Z}^m$ and $(L, \rho)$ be a partial character on
$\mathbb{Z}^m$. Then $I_{L, \rho}$ is complete intersection if and
only if for every face $\mathcal{F}$ of $\sigma_L$ the lattice
ideal $I_{L\cap \mathbb{Z}^{\mathbb{E}_{\mathcal{F}}},\rho}$ is
complete intersection.
\end{thm1}
\noindent \textbf{Proof.} ($\Leftarrow$) For
$\mathcal{F}=\sigma_L$ we have that $I_{L\cap
\mathbb{Z}^{\mathbb{E}_{\sigma_{L}}},\rho}=I_{L, \rho}$ is
complete
intersection.\\
($\Rightarrow$) Suppose that $I_{L, \rho}=(B({\bf
u}_1),\ldots,B({\bf u}_r))$ is complete intersection, where
$r=rank(L)$. Given a set $\{{\bf v}_1,\ldots,{\bf v}_q\} \subset
\mathbb{Z}^m$, we shall denote by $M({\bf v}_1,\ldots,{\bf v}_q)$
the $q \times m$ matrix whose rows are the vectors ${\bf
v}_1,\ldots,{\bf v}_q$. From Theorem 3.9 in \cite{MT} we have that
the matrix $M({\bf u}_1,\ldots,{\bf u}_r)$ is mixed dominating
while from Theorem 3.5 we deduce that $L \cap
\mathbb{Z}^{\mathbb{E}_{\mathcal{F}}}=\sum _{{\bf u}_i \in
\mathbb{Z}^{\mathbb{E}_{\mathcal{F}}}} \mathbb{Z}{\bf u}_i$. Set
$U=\{{\bf u}_1,\ldots,{\bf u}_r\}$ and $U \cap
\mathbb{Z}^{\mathbb{E}_{\mathcal{F}}}=\{{\bf u}_{i_1},\ldots,{\bf
u}_{i_k}\}$. Recall that a matrix $M$ with coefficients in
$\mathbb{Z}$ is called mixed dominating if it is mixed, i.e. every
row has a positive and negative entry, and also it does not
contain any square mixed submatrix. By Corollary 2.8 in \cite{FS}
the vectors ${\bf u}_1,\ldots,{\bf u}_r$ are linearly independent,
so in particular the set $U \cap
\mathbb{Z}^{\mathbb{E}_{\mathcal{F}}}$ is linearly independent.
Thus it is a $\mathbb{Z}$-base for the lattice $L \cap
\mathbb{Z}^{\mathbb{E}_{\mathcal{F}}}$. Also $M({\bf
u}_1,\ldots,{\bf u}_r)$ is mixed dominating, so $M({\bf
u}_{i_1},\ldots,{\bf u}_{i_k})$ is mixed dominating. Thus again
from Theorem 3.9 in \cite{MT} we have that $I_{L\cap
\mathbb{Z}^{\mathbb{E}_{\mathcal{F}}},\rho}$ is complete
intersection. \hfill $\square$\\

Next we will use the results of section 3 to generalize Theorem 1
and Theorem 2 of \cite{BMTh}.
\begin{prop1}Let $E_1\subset E_2\subset \{1,\ldots ,m\}$, $L \subset
\mathbb{Z}^m$ be a lattice and $L_1=L\cap \mathbb{Z}^{E_1}$,
$L_2=L\cap \mathbb{Z}^{E_2}$. Every $\mathbb{Z}$-basis (resp.
spanning set) of $L_1$ can be extended to a $\mathbb{Z}$-basis
(resp. spanning set) of $L_2$.
\end{prop1}
\noindent \textbf{Proof.} It is enough to consider the case that
$E_2=E_1\cup \{ i\}$. Let $\{{\bf u}_1,\ldots ,{\bf u}_r\}$ be a
$\mathbb{Z}$-basis of $L_1$. There are two cases.
\begin{enumerate} \item Every element ${\bf v} \in L_2$ has the $i$-coordinate equal to
zero. Then ${\rm supp}({\bf v}) \subset E_{1}$, so ${\bf v} \in
L\cap \mathbb{Z}^{E_1}=L_1$. Thus $L_2 \subset L_1$ and therefore
$L_2=L_1$. \item There are elements in $L_2$ with $i$-coordinate
different from zero. Choose ${\bf u}_{r+1} \in L_2$ such that its
$i$-th coordinate $({\bf u}_{r+1})_i$ is positive and this is the
least possible $i$-coordinate among all elements of $L_2$. Remark
that ${\bf u}_{r+1} \notin L_1$. If ${\bf u}_{r+1} \in L_{1}$,
then there are $\lambda_{1},\ldots,\lambda_{r} \in \mathbb{Z}$
such that ${\bf u}_{r+1}=\lambda_{1}{\bf
u}_1+\cdots+\lambda_{r}{\bf u}_r$. But $({\bf u}_{r+1})_{i}>0$
while $({\bf u}_1)_{i}=\cdots=({\bf u}_r)_{i}=0$. \newline Let
${\bf v}\in L_2$. If $({\bf v})_{i}=0$, then ${\bf v} \in L_1$.
Suppose that $({\bf v})_{i}>0$, the case that $({\bf v})_{i}<0$ is
similar. Divide $({\bf v})_{i}$ with $({\bf u}_{r+1})_{i}$, so
there exist $\lambda,\mu \in \mathbb{Z}$ such that $({\bf
v})_{i}=\lambda({\bf u}_{r+1})_{i}+ \mu $ and $0\leq \mu < ({\bf
u}_{r+1})_{i}$. Set ${\bf w}={\bf v}-\lambda{\bf u}_{r+1} \in
L_2$. Remark that $({\bf w})_{i}=\mu \geq 0$. Since ${\bf w} \in
L_{2}$ and $({\bf w})_{i} \geq 0$, we have, from the choice of
${\bf u}_{r+1}$, that $\mu =0$. Thus ${\bf w} \in L \cap
\mathbb{Z}^{E_1}=L_1$ and therefore ${\bf v} \in \sum_{i=1}^{r+1}
\mathbb{Z} {\bf u}_i$. It remains to prove that the vectors ${\bf
u}_1,\ldots,{\bf u}_r,{\bf u}_{r+1}$ are $\mathbb{Z}$-linearly
independent. Consider an equality
$$\lambda_{1}{\bf u}_1+\cdots+\lambda_{r}{\bf u}_r+\lambda_{r+1}{\bf u}_{r+1}={\mathbf
0}.$$Then $\lambda_{1}({\bf u}_1)_i+\cdots+\lambda_{r}({\bf
u}_r)_i+\lambda_{r+1}({\bf u}_{r+1})_i=0$ and therefore
$\lambda_{r+1}=0$. But ${\bf u}_1,\ldots,{\bf u}_r$ are
$\mathbb{Z}$-linearly independent, so
$\lambda_{1}=\cdots=\lambda_{r}=0$. \hfill $\square$
\end{enumerate}

\begin{def1} {\rm We say that a configuration of vectors $A=\{{\bf a}_1,\dots , {\bf a}_m\}$ is full if
\begin{enumerate} \item the cone $\sigma =pos_{\mathbb{Q}}(A)$ is an n-dimensional simplex
 cone generated by ${\bf a}_1,\dots , {\bf a}_n$, i.e. $\sigma=pos_{\mathbb{Q}}({\bf a}_{1},\ldots,{\bf a}_{n})$,
 and \item for each
 $n<i\leq m$ we have ${\mathcal{ F}}_{{\bf a}_i}=\sigma$. This means that the
 vectors ${\bf a}_i$, $n<i\leq m$, are in the relative
 interior $relint_{\mathbb{Q}}({\bf a}_{1},\ldots,{\bf a}_{n})$ of $\sigma$,
 which is the set of all strictly positive rational linear combinations of ${\bf a}_{1},\ldots,{\bf a}_{n}$.
 \end{enumerate}}
\end{def1}
\begin{thm1} Let $(L,\rho)$ be a partial character on $\mathbb{Z}^m$  such that
the configuration of vectors associated to the lattice $L$ is
full, then \begin{enumerate} \item in characteristic zero
$I_{L,\rho}$ is either a set theoretic complete intersection or an
almost set theoretic complete intersection \item in positive
characteristic $I_{L,\rho}$ is a set theoretic complete
intersection.
\end{enumerate}
\end{thm1}
\noindent \textbf{Proof.} Let $A= \{{\bf a}_{1},\ldots,{\bf
a}_{m}\} \subset \mathbb{Z}^{n}$ be the configuration of vectors
associated to $L$, then every face of $\sigma_{L}$ is of the form
$pos_{\mathbb{Q}}({\bf a}_{i}| i \in E)$, for a subset $E$ of
$\{1,\ldots,n\}$. Moreover for every face
$\mathcal{F}\not=\sigma_L $ we have that $L \cap
\mathbb{Z}^{\mathbb{E}_{\mathcal{F}}}=\{{\bf 0}\}$, since $A$ is
full. Applying Proposition 4.2 to the lattices $L \cap
\mathbb{Z}^{E_i}$, where $E_{i}=\{1,\ldots,i\}$ and $n<i\leq m$,
we obtain a $\mathbb{Z}$-basis $\{{\bf u}_{n+1}, \dots,{\bf
u}_{m}\}$ of the lattice $L$. Remark that ${\rm supp}({\bf u}_i)
\subset E_i$, the $i$-coordinate of ${\bf u}_i$ is positive and
also this is the least possible $i$-coordinate among all elements
of $L \cap \mathbb{Z}^{E_i}$. In addition, using the fact that for
every $i \in \{n+1,\ldots,m\}$ it holds ${\mathcal{F}}_{{\bf
a}_i}=\sigma_{L}$, we take binomials $B({\bf v}_i) \in
I_{L,\rho}$, $n<i\leq m$, with ${\rm supp}(({\bf v}_{i})_+)=\{1,
\dots, n\}$ and ${\rm supp}(({\bf v}_{i})_-)=\{i\}$. Note that
${\bf u}_{n+1}$ is in that form, so we will consider ${\bf
v}_{n+1}:={\bf u}_{n+1}$. Now we will distinguish two cases for
the characteristic of the field $K$.
\begin{enumerate} \item The characteristic of $K$ is equal to zero. We will
construct $m-n+1$ binomials that generate the radical of
$I_{L,\rho}$ up to radical. Given an index $i \in
\{n+3,\ldots,m\}$, there are appropriate large positive integers
$r_{n+1},\ldots,r_{i-2}$ such that for every $j \in
\{n+1,\ldots,i-2\}$ the $j$-coordinate of ${\bf w}_i={\bf
u}_i-\sum_{j=n+1}^{i-2}r_j{\bf v}_j$ is positive. Furthermore,
there exist a big enough positive integer $r$ such that for every
$j \in \{1,\ldots,n\}$ the $j$-coordinate of ${\bf z}_i={\bf
w}_{i}+r{\bf v}_{i-1}$ is positive, while also the $j$-coordinate,
$n+1 \leq j \leq i-2$, of ${\bf z}_i$ is positive, since are the
same as ${\bf w}_i$. Notice that the vectors ${\bf z}_i \in L$ and
${\bf u}_{i}$ have the same $i$-coordinate and also ${\rm
supp}({\bf z}_i) = E_i$. But the lattice $L$ is positive, so $L
\cap \mathbb{N}^{m}=\{{\bf 0}\}$ and therefore the
$(i-1)$-coordinate of ${\bf z}_i$ is negative. Set ${\bf
z}_{n+1}={\bf u}_{n+1}$ and ${\bf z}_{n+2}={\bf u}_{n+2}$. From
the proof of Proposition 4.2 we have that $\{{\bf
z}_{n+1},\ldots,{\bf z}_{m}\}$ is a $\mathbb{Z}$-basis of $L$. Let
${\bf z}_{m+1}={\bf v}_m$, then for every $i \in \{n+1,
n+3,\ldots,m+1\}$ the binomial $B({\bf z}_i)$ is of the form
$B({\bf z}_i)=x_{i-1}^{g_{i-1}}-N_{i-1}$ where $N_{i-1}$ is a
monomial not containing the variable $x_{i-1}$. Thus the set
$\{B({\bf z}_{n+1}), B({\bf z}_{n+3}),\ldots,B({\bf z}_{m+1})\}$
is a cover of $A$. Using Theorem 3.5 we take that the $m-n+1$
binomials $B({\bf z}_i)$, where $n<i\leq m+1$, generate the
radical of the lattice ideal up to radical. \item The
characteristic of $K$ is equal to $p>0$. We are going to construct
$m-n$ binomials that generate the radical of $I_{L,\rho}$ up to
radical. Given an index $i \in \{n+2,\ldots,m\}$, there are
appropriate large positive integers $r_{n+1},\ldots,r_{i-1}$ such
that for every $j \in \{n+1,\ldots,i-1\}$ the $j$-coordinate of
${\bf w}_i={\bf u}_i+\sum_{j=n+1}^{i-1}r_j{\bf v}_j$ is negative.
Note that
\begin{enumerate} \item the $i$-coordinate of ${\bf w}_i$
coincides with the $i$-coordinate of ${\bf u}_i$ and both of them
are positive, \item the $i$-coordinate of ${\bf v}_i$ is equal to
$-t({\bf u}_i)_i$ for a positive integer $t$, where $({\bf
u}_i)_{i}$ is the $i$-coordinate of ${\bf u}_{i}$, since $({\bf
u}_i)_{i}$ is the least possible $i$-coordinate among all elements
of $L \cap \mathbb{Z}^{E_{i}}$.
\end{enumerate}
Given a power $p^{k}$ of $p$, there exist integers $r,s \in
\mathbb{Z}$ such that $p^{k}=rt+s$ and $0 \leq s<t$. We choose a
big enough $k$ such that for every $j \in \{1,\ldots,n\}$ the
$j$-coordinate of ${\bf y}_i:=s{\bf w}_{i}-r{\bf v}_{i} \in L$,
$n+1<i \leq m$, is negative. Therefore the binomial $B({\bf y}_i)$
is of the form $B({\bf y}_i)=x_{i}^{g_{i}}-N_{i}$, where $N_{i}$
is a monomial not containing the variable $x_{i}$, ${\rm
supp}({\bf y}_i)=E_i$ and the $i$-coordinate of ${\bf y}_i$ is
$p^k({\bf u}_i)_i$. Let ${\bf y}_{n+1}={\bf u}_{n+1}$, then the
set $\{B({\bf y}_{n+1}),\ldots,B({\bf y}_{m})\}$ is a cover of
$A$. It remains to prove that
$$(L\cap
\mathbb{Z}^{E_i}):p^{\infty}=(\sum _{j=n+1}^{i}\mathbb{Z}{\bf
y}_j):p^{\infty}, \ \textrm{for every} \ i=n+1,\ldots,m.$$ The
proof is obvious for $i=n+1$, since $$L\cap
\mathbb{Z}^{E_{n+1}}=\mathbb{Z}{\bf u}_{n+1}=\mathbb{Z}{\bf
y}_{n+1}.$$ Assume that $(L\cap
\mathbb{Z}^{E_{i-1}}):p^{\infty}=(\sum
_{j=n+1}^{i-1}\mathbb{Z}{\bf y}_j):p^{\infty}$. We have $$(\sum
_{j=n+1}^{i}\mathbb{Z}{\bf y}_j):p^{\infty} \subset (L\cap
\mathbb{Z}^{E_i}):p^{\infty},$$ since $(\sum
_{j=n+1}^{i}\mathbb{Z}{\bf y}_j) \subset L\cap \mathbb{Z}^{E_i}$.
Now we will prove that $$(L\cap \mathbb{Z}^{E_i}):p^{\infty}
\subset (\sum _{j=n+1}^{i}\mathbb{Z}{\bf y}_j):p^{\infty}.$$ Let
$g \in (L\cap \mathbb{Z}^{E_i}):p^{\infty}$, then $p^{b}g \in
L\cap \mathbb{Z}^{E_i}$, for a $b \in \mathbb{N}$, and therefore
$p^{b}g=\lambda_{n+1}{\bf u}_{n+1}+\cdots+\lambda_{i}{\bf u}_{i}$
for some integers $\lambda_{n+1},\ldots,\lambda_{i}$. Remark that
$\lambda_{n+1}{\bf u}_{n+1}+\cdots+\lambda_{i-1}{\bf u}_{i-1}$
belongs to $L \cap \mathbb{Z}^{E_{i-1}}$. This means that there is
a $c \in \mathbb{N}$ such that $p^{c}(\lambda_{n+1}{\bf
u}_{n+1}+\cdots+\lambda_{i-1}{\bf u}_{i-1})$ belongs to $(\sum
_{j=n+1}^{i-1}\mathbb{Z}{\bf y}_j):p^{\infty}$. Moreover ${\bf
y}_i-p^k{\bf u}_i$ belongs to $L\cap \mathbb{Z}^{E_{i-1}}$, so
there is a natural number $d$ such that $p^{d}({\bf y}_{i}-p^k{\bf
u}_i) \in \sum _{j=n+1}^{i-1}\mathbb{Z}{\bf y}_j$ and therefore
$p^{d+k}{\bf u}_i \in \sum _{j=n+1}^{i}\mathbb{Z}{\bf y}_j$. Thus
$p^{c+d+k}g$ belongs to $\sum _{j=n+1}^{i}\mathbb{Z}{\bf y}_j$.
Using Theorem 3.5 we take that the $m-n$ binomials $B({\bf y}_i)$,
where $n<i\leq m$, generate the radical of $I_{L, \rho}$ up to
radical. \hfill $\square$
\end{enumerate}

\begin{ex1} {\rm Let $m$ be a fixed integer number greater than or equal
to $8$ and let $L(m)$ be the sublattice of $\mathbb{Z}^4$
generated by ${\bf e}_1=(m+2q-3,-m+2q+5,-1,-1)$, ${\bf
e}_2=(-m-2q+5,m-2q-3,-1,-1)$, ${\bf e}_3=(-m-2q+5,-1,m-3,-1)$,
where $q=0$ when $m$ is even and $q=1$ otherwise. It is easy to
check that $L(m)$ has rank 3; thus, $I_{L(m), 1}$ is a lattice
ideal in $K[x_{1},\ldots,x_{4}]$ of codimension $3$ where $K$ is a
field of any characteristic.
 Moreover $L(m)$ is not saturated because $(2,2,-2,-2)={\bf e}_1+{\bf e}_2 \in L(m)$ and
 it is easy to check that  $(1,1,-1,-1) \notin L(m)$; therefore, the lattice ideal $I_{L(m), 1}$ is never
 toric. In \cite{Oj} Ojeda proved that the lattice ideal $I_{L(m), 1}$ is generic and minimally generated  by $m$ elements, i.e. $$f_{1}=x_{1}^{m+2q-3}-x_{2}^{m-2q-5}x_{3}x_{4}, \ \ f_{2}=x_{2}^{m-2q-3}-x_{1}^{m+2q-5}x_{3}x_{4},$$ $$f_{3}=x_{3}^{m-3}-x_{1}^{m+2q-5}x_{2}x_{4}$$ $$f_{12}=x_{1}^2x_{2}^2-x_{3}^2x_{4}^2$$ $$f_{23}^{i}=x_{2}^{2i-1}x_{3}^{m-2i-3}-x_{1}^{m+2q-2i-5}x_{4}^{2i+1}, \ i=1,\ldots,\frac{m}{2}+\frac{q}{2}-3$$
$$f_{13}^{j}=x_{1}^{2j}x_{3}^{m-2j-2}-x_{2}^{m-2q-2j-2}x_{4}^{2j}, \ j=1,\ldots,\frac{m}{2}-\frac{q}{2}-2$$ and $$f_{123}=
\left\{ \begin{array}{ll}
         x_{4}^{m-3}-x_{1}x_{2}^{m-5}x_{3} & \mbox{if m is even}\\
        x_{4}^{m-3}-x_{1}^{m-3}x_{2}x_{3} & \mbox{if m is odd}.\end{array} \right.$$
On the other hand for every integer $m \geq 8$ consider the
homomorphism $\phi_{m}:\mathbb{Z}^4 \rightarrow \mathbb{Z}$
defined by $\phi_{m}(a,b,c,d)=a+b+c+d$ if $m$ is even, and
$\phi_{m}(a,b,c,d)=a(m-6)+b(m-2)+c(m-4)+d(m-4)$ if $m$ is odd. By
direct computation we can check that $L(m)\subset\ker \phi_{m}$,
so that $Sat(L(m))=\ker \phi_{m}$. We can associate to $L(m)$ the
rational polyhedral cone $\sigma_{L(m)}=pos_{\mathbb{Q}}(A)$ where
\begin{enumerate} \item $A=\{\,1,1,1,1\}$ when $m$ is even and
\item $A=\{\,m-6,m-2,m-4,m-4\}$ when $m$ is odd.
\end{enumerate}
In both cases the cone $\sigma_{L(m)}=pos_{\mathbb{Q}}(A)$ has
only one nonzero  face $\mathcal{F}=pos_{\mathbb{Q}}(A)$. So
${L(m)} \cap
(\mathbb{Z}^4)^{\mathbb{E}_{\sigma_{{L(m)}}}}={L(m)}=<{\bf
e}_1,{\bf e}_2,{\bf e}_3>$. To cover $A$ we need also the binomial
$f_{123}$. So, from Theorem 3.5, we have that
$rad(I_{L(m),1})=rad(f_{1},f_{2},f_{3},f_{123})$. Therefore
$I_{L(m),1}$ is an almost set theoretic complete intersection in
the characteristic zero case. In positive characteristic we have,
from Theorem 4.4, that it is set theoretic complete intersection
and the $3$ binomials which generate $rad(I_{L(m),1})$ up to
radical depend on the characteristic.}
\end{ex1}

\end{document}